\documentclass[12pt,a4paper]{amsart}

\usepackage{amsmath, amssymb, amsthm, xcolor}
\usepackage{hyperref}

\newtheorem{definition}{Definition}[section]
\newtheorem{theorem}{Theorem}[section]
\newtheorem{corollary}[theorem]{Corollary}

\newtheorem{proposition}{Proposition}[section]

\newcommand{\SubX}{Sub_X} 
\newcommand{\That}{\hat{T}}
\newcommand{\VOpen}[1]{\langle #1 \rangle}

\begin{document}

\begin{center}
    
\large{\textbf{$k$-type Chaos for Induced Group Actions on Hyperspaces}}

\small{Anshid Aboobacker}

\textit{Department of Mathematics, BITS Pilani Hyderabad Campus, India}

anshidaboobackerk@gmail.com
\end{center}

\noindent{\bf Abstract.}
This paper investigates the correlation between $k$-type dynamical properties of $\mathbb{Z}^d$-actions on compact metric spaces and their induced actions on the corresponding hyperspaces. We extend the classical results from discrete dynamical systems and general group actions to the specific setting of $k$-type dynamics. Specifically, we define and study $k$-type transitivity, $k$-type mixing, $k$-type weak mixing, and $k$-type Li-Yorke chaos for induced hyperspace actions, establishing that these properties transfer from the base system to the hyperspace under appropriate conditions.

\vspace{0.5cm}

\section{Introduction}

Given a topological space $X$, its hyperspace $\mathcal{K}(X)$ is defined as the collection of all non-empty compact subsets of $X$. The study of dynamical systems on hyperspaces has a rich history, offering insights into the ``collective'' behavior of a system when individual points are replaced by subsets. Traditionally, much of this work has focused on discrete dynamical systems generated by a single continuous map $f: X \to X$, where the induced map $\hat{f}$ acts on the hyperspace $\mathcal{K}(X)$ (see, e.g., \cite{Akin2017, Banks2005, prajapati2020expansive, RomanFlores2003, Sharma2010}). These studies have revealed profound connections, such as the equivalence between certain chaotic properties of $f$ and those of $\hat{f}$. Recently, there have been several attempts to extend these ideas to more general settings involving group actions and their induced actions on hyperspaces (see, for example, \cite{huang2022topological}).

In this paper, we focus on dynamical systems generated by $\mathbb{Z}^d$-actions. Throughout, $(X,d)$ denotes a compact metric space, and we consider continuous actions of $\mathbb{Z}^d$ on $X$. We define the induced action on the hyperspace $\SubX$, the collection of all closed subsets of $X$. Since $X$ is compact and metric, we have $\SubX = \mathcal{K}(X)$; however, we retain the notation $\SubX$ to emphasize the hyperspace viewpoint. Our main objective is to investigate how the dynamical behavior of the $\mathbb{Z}^d$-action on $X$ is reflected in, and influenced by, the corresponding induced action on $\SubX$.

A fundamental difficulty in passing from single-map dynamics to $\mathbb{Z}^d$-actions lies in the absence of a natural linear notion of time. For a single map, iterates $f^n$ and asymptotic behavior as $n \to \infty$ are well-defined. In contrast, for $\mathbb{Z}^d$-actions, dynamics evolve along multiple directions, and the idea of ``large iterates'' must be interpreted different. To address this issue, we adopt the notion of $k$-type orders on $\mathbb{Z}^d$ introduced by Oprocha \cite{oprocha2007chain} and further developed by Shah and Das \cite{shah2015different, shah2015note}. 

A $\mathbb{Z}^d$-action on a topological space $X$ is defined as a continuous mapping
$T : \mathbb{Z}^d \times X \to X$
satisfying the identities $T^{0}(x)=x$ and
$T^{n+m}(x)=T^{n}\bigl(T^{m}(x)\bigr)$
for all $x\in X$ and $n,m\in\mathbb{Z}^d$, where $T^{n}(x)$ is a shorthand for $T(n,x)$.
Let $T:\mathbb{Z}^d \times X \to X$ be such an action and fix an integer $k \in \{1,2,\dots,2^d\}$. Denote by $k^b$ the binary representation of $k-1$ written in $d$ digits, that is,
$k-1 = \sum_{i=1}^{d} k_i^b\, 2^{\,i-1}$,
where $k^b=(k_1^b,\dots,k_d^b)\in\{0,1\}^d$.
For vectors $x,y\in\mathbb{Z}^d$, we define the relation $x >^{k} y$ if
$(-1)^{k_i^b} x_i > (-1)^{k_i^b} y_i$ for all $i=1,\dots,d$,
where $x=(x_1,\dots,x_d)$ and $y=(y_1,\dots,y_d)$.

Using this $k$-type order, we study $k$-type transitivity, $k$-type weak mixing, $k$-type mixing, $k$-type poximal pairs, $k$-type asymptotic pairs, and $k$-type Li-Yorke chaos as they appear for $\mathbb{Z}^d$-actions in the literature \cite{aboobacker2025ktype,oprocha2007chain,shah2015different,shah2015note}. Our focus is on understanding how these $k$-type dynamical properties of a $\mathbb{Z}^d$-action on $X$ relate to the corresponding properties of the induced action on the hyperspace $\SubX$.

The organization of the paper is as follows. In Section~\ref{sec:preliminaries}, we collect the necessary preliminaries on hyperspace topology and the $k$-type dynamical notions used throughout the paper. Section~\ref{sec:mainresults} contains our main results, where we establish implications and equivalences between $k$-type dynamical properties of $(X,T)$ and those of the induced action on $(\SubX,\hat{T})$.


\section{Preliminaries}
\label{sec:preliminaries}

\subsection{The Hyperspace $\SubX$}

The hyperspace $\SubX$ is endowed with the \emph{Hausdorff metric} $H$, defined for $A, B \in \SubX$ as:
\[
H(A, B) = \max \left\{ \sup_{a \in A} d(a, B), \sup_{b \in B} d(b, A) \right\},
\]
where $d(x, S) = \inf\limits_{s \in S} d(x, s)$ for $x \in X$ and $S \subseteq X$.
It is well known that $(\SubX, H)$ is also a compact metric space since $(X,d)$ is a compact metric space.

The \emph{Vietoris topology} on $\SubX$ is generated by a basis of sets of the form:
\[
\VOpen{U_1, \dots, U_n} := \left\{ A \in \SubX : A \subseteq \bigcup_{i=1}^n U_i \text{ and } A \cap U_i \neq \emptyset \text{ for all } 0 \le i \le n \right\},
\]
where $U_1, \dots, U_n$ are non-empty open subsets of $X$. Since $X$ is a compact metric space, the Vietoris topology coincides with the topology induced by the Hausdorff metric on $\SubX$. Throughout this paper, $\SubX$ is endowed with the Hausdorff metric.

Note that when considering two basic open sets $\mathcal{U}$ and $\mathcal{V}$ in $\SubX$, we may assume without loss of generality that they are each generated by the same number of open subsets of $X$. That is, if
\[
\mathcal{U}=\VOpen{U_1,\dots,U_m},\qquad
\mathcal{V}=\VOpen{V_1,\dots,V_n},
\]
then setting $k=\max\{m,n\}$ and extending the smaller family by putting $U_i=X$ or $V_j=X$ for the appropriate indices yields representations of $\mathcal{U}$ and $\mathcal{V}$ as
\[
\mathcal{U}=\VOpen{U_1,\dots,U_k},\qquad
\mathcal{V}=\VOpen{V_1,\dots,V_k}.
\]

\subsection{Induced $\mathbb{Z}^d$-Action on $\SubX$}

Let $T: \mathbb{Z}^d \times X \to X$ be a $\mathbb{Z}^d$-action. This naturally induces an action on $\SubX$. For each $m \in \mathbb{Z}^d$, we define the induced map $\hat{T}^m: \SubX \to \SubX$ by:
\[
\hat{T}^m(A) = T^m(A) = \{T^m(x) : x \in A\} \quad \text{for all } A \in \SubX.
\]
Since $T^m$ is a homeomorphism and $A$ is compact, $T^m(A)$ is compact. The map $\hat{T}^m$ is a homeomorphism on $\SubX$, and the group property $\hat{T}^{m+n} = \hat{T}^m \circ \hat{T}^n$ holds. Thus, $(\SubX, \hat{T})$ is a $\mathbb{Z}^d$-action.

\subsection{$k$-type Dynamical Notions}

We now define some dynamical properties for general group actions that are needed for this article. It may be noted that all these are available in the literature.  Let $k \in \{1, \dots, 2^d\}$. 

\begin{definition}
	The action $T$ is \emph{$k$-type transitive} if for any pair of non-empty open sets $U, V \subseteq X$, there exists $m \in \mathbb{Z}^d$ with $m \ge^k 0$ such that $T^m(U) \cap V \neq \emptyset$.
\end{definition}

\begin{definition}
	The action $T$ is \emph{$k$-type weakly mixing} if the product action $T \times T$ on $X \times X$ is $k$-type transitive. Equivalently, for any four non-empty open sets $U_1, U_2, V_1, V_2 \subseteq X$, there exists $m \ge^k 0$ such that
	\[
	T^m(U_1) \cap V_1 \neq \emptyset \quad \text{and} \quad T^m(U_2) \cap V_2 \neq \emptyset.
	\]
\end{definition}

\begin{definition}
	The action $T$ is \emph{$k$-type mixing} if for any pair of non-empty open sets $U, V \subseteq X$, there exists a finite subset $F \subset \{n \in \mathbb{Z}^d : n \ge^k 0\}$ such that 
    \[
    T^m(U) \cap V \neq \emptyset
    \]
    for all $m \ge^k 0$ with $m \notin F$.
\end{definition}





For an $\epsilon>0$, define 
\[ 
V_{\epsilon} = \{ (x,y)\in X\times X : d(x,y) < \epsilon \}, \]
\[
\overline{V_{\epsilon}} = \{ (x,y)\in X\times X : d(x,y) \leq \epsilon \}. \]

\begin{definition} Let $(X,T)$ be a dynamical system given by a $\mathbb{Z}^d$-action $T$ on $X$. The set of all \emph{$k-$type proximal pairs} of $T$ is defined as

\[k-Prox(T) = \bigcap_{\epsilon>0}\bigcup_{\substack{{n}>^k {0} \\ {n}\in \mathbb{Z}^d}} \left[ T^{{-n}}\times T^{{-n}} (V_{\epsilon}) \right]. \]
Any pair $(x,y) \in k-Prox(T) $ is called a $k-$type proximal pair of $T$.
\end{definition} 

$k-$type proximal pairs of $T$ has an equivalent definition as given in \cite{aboobacker2025ktype}.

\begin{proposition}
    $(x,y) \in k-Prox(T)$ if and only if there exists a sequence $({t_s})_{s\in \mathbb{N}} $ in $\mathbb{Z}^d$ with $ {t_{s+1}} >^k {t_{s}}  $ such that $d(T^{{t_s}}(x) , T^{{t_s}}(y)) \rightarrow 0 $ as $s\rightarrow \infty$.
\end{proposition}

\begin{definition} For any $\epsilon>0$, \[ k-Asym_{\epsilon}(T)= \bigcup_{\substack{{n}>^k {0} \\ {n}\in \mathbb{Z}^d}} T^{{-n}}\times T^{{-n}} \left( \bigcap_{\substack{{m}>^k {0} \\ {m}\in \mathbb{Z}^d}} [ T^{{-m}}\times T^{{-m}} (\overline{V_{\epsilon}})] \right) \]

and  \[ k-Asym(T) = \bigcap_{\epsilon>0} k-Asym_{\epsilon}(T).   \]
We call elements $(x,y)$ in $k-Asym(T)$ as \emph{$k-$type asymptotic pairs} of $T$.

\end{definition}

\begin{definition}
    Two points $x,y\in X$ are said to form a \emph{$k-$type Li-Yorke pair} if $(x,y)\in k-Prox(T)\setminus k-Asym(T)$ and the set of all $k$-type Li-Yorke pairs is denoted by $k-LY(T)$. 
    A subset $S$ of $X$ is said to be a \emph{$k-$ scrambled set} if any two distinct pairs of points in $S$ form a $k-$type Li-Yorke pair.
    A dynamical system $(X,T)$ is said to be \emph{$k$-type Li-Yorke chaotic} if $X$ contains an uncountable scrambled set.
\end{definition}

\medskip

The corresponding $k$-type dynamical properties for the induced action $(\SubX, \hat{T})$ are defined analogously by replacing $X$ with $\SubX$, $T$ with $\hat{T}$, points $x$ with compact sets $A$, and the metric $d$ with the Hausdorff metric $H$.

\section{Main Results}
\label{sec:mainresults}

In this section, we establish fundamental relationships between the $k$-type dynamical properties of a $\mathbb{Z}^d$-action $(X,T)$ and those of its induced action $(\SubX,\hat{T})$. Many of the results presented here are adaptations of known theorems for $\mathbb{Z}$-actions, suitably modified to account for the directional nature of $\mathbb{Z}^d$-actions.

\subsection{$k$-type Transitivity}

We begin by examining the relationship between $k$-type transitivity of the action on $X$ and $k$-type transitivity of the induced action on the hyperspace $\SubX$.

Following the approach of Rom\'an-Flores \cite{RomanFlores2003}, for any open set $U \subseteq X$, we define its \emph{extension} to the hyperspace $\SubX$ as the collection of all non-empty compact subsets of $X$ that are entirely contained in $U$. We will also use Proposition~\ref{prop:e_is_open}, which is proved in \cite{RomanFlores2003}.

\begin{definition}
For any subset $S \subseteq X$, its extension to the hyperspace $\SubX$ is defined by
\[
e(S) := \{K \in \SubX \mid K \subseteq S\}.
\]
\end{definition}

\noindent The following proposition shows that these extended sets are open in the Vietoris topology when the original set is open, and describes how they behave under intersections.

\begin{proposition}
\label{prop:e_is_open}
If $U$ is a non-empty open subset of $X$, then $e(U)$ is a non-empty open subset of $\SubX$. Moreover, for any subsets $U,V \subseteq X$,
\[
e(U) \cap e(V) = e(U \cap V).
\]
\end{proposition}

\noindent The next proposition describes how extended sets interact with the induced $\mathbb{Z}^d$-action.

\begin{proposition}
\label{prop:e_properties}
Let $U \subseteq X$ and let $m \in \mathbb{Z}^d$. Then
\[
\hat{T}^{m}\bigl(e(U)\bigr) \subseteq e\bigl(T^{m}(U)\bigr).
\]
\end{proposition}

\begin{proof}
Let $A \in \hat{T}^{m}(e(U))$. By definition of the induced action, there exists a set $K \in e(U)$ such that
\[
A = \hat{T}^{m}(K) = T^{m}(K).
\]
Since $K \subseteq U$ and $T^{m}$ is a homeomorphism, it follows that
\[
T^{m}(K) \subseteq T^{m}(U).
\]
Hence $A \subseteq T^{m}(U)$, which implies $A \in e(T^{m}(U))$. This proves the inclusion.
\end{proof}

We now prove the main result relating $k$-type transitivity of the induced action to $k$-type transitivity of the base action.

\begin{theorem}
\label{thm:induced_k_trans_implies_base_k_trans}
If the induced action $(\SubX,\hat{T})$ is $k$-type transitive, then the $\mathbb{Z}^d$-action $(X,T)$ is $k$-type transitive.
\end{theorem}

\begin{proof}
Assume that $(\SubX,\hat{T})$ is $k$-type transitive. Let $U,V \subseteq X$ be arbitrary non-empty open sets. We must show that there exists an element $m \in \mathbb{Z}^d$ with $m \ge^k 0$ such that
\[
T^{m}(U) \cap V \neq \emptyset.
\]

By Proposition~\ref{prop:e_is_open}, the sets $e(U)$ and $e(V)$ are non-empty open subsets of $\SubX$. Since $(\SubX,\hat{T})$ is $k$-type transitive, there exists $m \in \mathbb{Z}^d$ with $m \ge^k 0$ such that
\[
\hat{T}^{m}\bigl(e(U)\bigr) \cap e(V) \neq \emptyset.
\]

By Proposition~\ref{prop:e_properties}, we have
\[
\hat{T}^{m}\bigl(e(U)\bigr) \subseteq e\bigl(T^{m}(U)\bigr),
\]
and hence
\[
e\bigl(T^{m}(U)\bigr) \cap e(V) \neq \emptyset.
\]
Using Proposition~\ref{prop:e_is_open}, this is equivalent to
\[
e\bigl(T^{m}(U) \cap V\bigr) \neq \emptyset.
\]

By definition, $e(S)$ is non-empty if and only if $S$ contains a non-empty compact subset, which implies that $S$ itself is non-empty. Therefore,
\[
T^{m}(U) \cap V \neq \emptyset.
\]
Since $U$ and $V$ were arbitrary non-empty open subsets of $X$, this shows that $(X,T)$ is $k$-type transitive.
\end{proof}

\subsection{$k$-type weakly mixing}
\noindent We begin by recording a property of $k$-type weakly mixing $\mathbb{Z}^d$-actions that will be used repeatedly in the sequel.

\begin{proposition}
\label{prop:k_weakmixing_iff_product}
Let $(X,T)$ be $k$-type weakly mixing and let 
$U_1, U_2, \dots, U_n,  V_1, V_2, \dots, V_n$
be non-empty open subsets of $X$. Then there exists an element 
$m \in \mathbb{Z}^d$ with $m \ge^k 0$ such that
\[
T^{m}(U_i) \cap V_i \neq \emptyset 
\quad \text{for all } i \in \{1,2,\dots,n\}.
\]
\end{proposition}

\begin{proof}
For non-empty open sets $U,V \subseteq X$, define
\[
N(U,V) := \{\, m \in \mathbb{Z}^d : m \ge^k 0 \text{ and } T^{m}(U) \cap V \neq \emptyset \,\}.
\]
We need to show that
\[
\bigcap_{i=1}^{n} N(U_i,V_i) \neq \emptyset.
\]

Since $(X,T)$ is $k$-type weakly mixing, each set $N(U,V)$ is non-empty.
We first show that there exist non-empty open sets $E$ and $F$ such that
\[
N(E,F) \subseteq N(U_1,V_1) \cap N(U_2,V_2).
\]

By $k$-type weak mixing, there exists $m \ge^k 0$ such that
\[
T^{m}(U_1) \cap U_2 \neq \emptyset
\quad \text{and} \quad
T^{m}(V_1) \cap V_2 \neq \emptyset.
\]
Define
\[
E := U_1 \cap T^{-m}(U_2),
\qquad
F := V_1 \cap T^{-m}(V_2).
\]
Then $E$ and $F$ are non-empty open subsets of $X$.

Let $m' \in N(E,F)$. Then
\[
T^{m'}(E) \cap F \neq \emptyset,
\]
that is,
\[
T^{m'}(U_1) \cap T^{m'-m}(U_2) \cap V_1 \cap T^{-m}(V_2) \neq \emptyset.
\]
This implies $T^{m'}(U_1) \cap V_1 \neq \emptyset$, hence $m' \in N(U_1,V_1)$.
Applying $T^{m}$ to the above intersection, we obtain
\[
T^{m'}(U_2) \cap V_2 \neq \emptyset,
\]
and therefore $m' \in N(U_2,V_2)$.
Thus,
\[
N(E,F) \subseteq N(U_1,V_1) \cap N(U_2,V_2).
\]

\noindent We now proceed inductively.
Suppose that for some $k_0 < n$, there exist non-empty open sets
$E_{k_0-1}$ and $F_{k_0-1}$ such that
\[
N(E_{k_0-1},F_{k_0-1})
\subseteq \bigcap_{i=1}^{k_0} N(U_i,V_i).
\]
Using the same argument as above, we can construct non-empty open sets
$E_{k_0}$ and $F_{k_0}$ satisfying
\[
N(E_{k_0},F_{k_0})
\subseteq N(E_{k_0-1},F_{k_0-1}) \cap N(U_{k_0+1},V_{k_0+1})
\subseteq \bigcap_{i=1}^{k_0+1} N(U_i,V_i).
\]
By induction, we obtain non-empty open sets $E_{n-1}$ and $F_{n-1}$ such that
\[
N(E_{n-1},F_{n-1})
\subseteq \bigcap_{i=1}^{n} N(U_i,V_i),
\]
and hence the intersection is non-empty.
\end{proof}

\begin{theorem}
\label{thm:k_main_equivalence}
The $\mathbb{Z}^d$-action $(X,T)$ is $k$-type weakly mixing if and only if the induced action $(\SubX,\That)$ is $k$-type weakly mixing.
\end{theorem}

\begin{proof}
Assume that $(X,T)$ is $k$-type weakly mixing. We show that the induced action $(\SubX,\That)$ is also $k$-type weakly mixing.

Let $\mathcal{U}_1, \mathcal{U}_2, \mathcal{V}_1, \mathcal{V}_2$ be non-empty open subsets of $\SubX$.
Without loss of generality, assume
\begin{align*}
\mathcal{U}_1 &= \VOpen{U_{1,1},\dots,U_{1,n}},
& \mathcal{V}_1 &= \VOpen{V_{1,1},\dots,V_{1,n}}, \\
\mathcal{U}_2 &= \VOpen{U_{2,1},\dots,U_{2,n}},
& \mathcal{V}_2 &= \VOpen{V_{2,1},\dots,V_{2,n}},
\end{align*}
where all $U_{i,j}$ and $V_{i,j}$ are non-empty open subsets of $X$.

By Proposition~\ref{prop:k_weakmixing_iff_product}, there exists $m \ge^k 0$ such that
\begin{enumerate}
\item[(a)] $T^{m}(U_{1,i}) \cap V_{1,i} \neq \emptyset$ for all $i$,
\item[(b)] $T^{m}(U_{2,i}) \cap V_{2,i} \neq \emptyset$ for all $i$.
\end{enumerate}

Define
\[
G := \VOpen{U_{1,1} \cap T^{-m}(V_{1,1}), \dots, U_{1,n} \cap T^{-m}(V_{1,n})},
\]
\[
H := \VOpen{U_{2,1} \cap T^{-m}(V_{2,1}), \dots, U_{2,n} \cap T^{-m}(V_{2,n})}.
\]
Then $G$ and $H$ are non-empty open subsets of $\SubX$.

If $A \in G$, then $A \in \mathcal{U}_1$ and $\That^{m}(A) \in \mathcal{V}_1$.
Similarly, if $B \in H$, then $B \in \mathcal{U}_2$ and $\That^{m}(B) \in \mathcal{V}_2$.
Thus,
\[
\That^{m}(\mathcal{U}_1) \cap \mathcal{V}_1 \neq \emptyset
\quad \text{and} \quad
\That^{m}(\mathcal{U}_2) \cap \mathcal{V}_2 \neq \emptyset,
\]
showing that $(\SubX,\That)$ is $k$-type weakly mixing.

\medskip

Conversely, assume that $(\SubX,\That)$ is $k$-type weakly mixing.
Let $U_1,V_1,U_2,V_2 \subseteq X$ be non-empty open sets, and define
\[
\mathcal{U}_1 := \VOpen{U_1}, \quad
\mathcal{U}_2 := \VOpen{U_2}, \quad
\mathcal{V}_1 := \VOpen{V_1}, \quad
\mathcal{V}_2 := \VOpen{V_2}.
\]
Then $\mathcal{U}_1, \mathcal{U}_2, \mathcal{V}_1, \mathcal{V}_2$ are non-empty open subsets of $\SubX$.
By $k$-type weak mixing of $(\SubX,\That)$, there exists $m \ge^k 0$ such that
\[
\That^{m}(\mathcal{U}_1) \cap \mathcal{V}_1 \neq \emptyset
\quad \text{and} \quad
\That^{m}(\mathcal{U}_2) \cap \mathcal{V}_2 \neq \emptyset.
\]
This implies
\[
T^{m}(U_1) \cap V_1 \neq \emptyset
\quad \text{and} \quad
T^{m}(U_2) \cap V_2 \neq \emptyset,
\]
and hence $(X,T)$ is $k$-type weakly mixing.
\end{proof}

\noindent Since a $k$-type weakly mixing action is $k$-type transitive, we obtain the following immediate corollaries.

\begin{corollary}
\label{cor:k_weak_implies_k_trans_SubX}
If the $\mathbb{Z}^d$-action $(X,T)$ is $k$-type weakly mixing, then the induced action $(\SubX,\That)$ is $k$-type transitive.
\end{corollary}

\begin{corollary}
\label{cor:k_weak_implies_k_trans_X}
If the induced action $(\SubX,\That)$ is $k$-type weakly mixing, then the $\mathbb{Z}^d$-action $(X,T)$ is $k$-type transitive.
\end{corollary}

\begin{theorem}
\label{thm:k_mixing_equivalence}
The $\mathbb{Z}^d$-action $(X,T)$ is $k$-type mixing if and only if the induced action $(\SubX,\That)$ is $k$-type mixing.
\end{theorem}

\begin{proof}
Assume that the action $(X,T)$ is $k$-type mixing.
Let $\mathcal{U}, \mathcal{V}$ be non-empty open subsets of $\SubX$. Suppose
\[
\mathcal{U} = \VOpen{U_1,\dots,U_n}
\quad \text{and} \quad
\mathcal{V} = \VOpen{V_1,\dots,V_n},
\]
where each $U_i$ and $V_i$ is a non-empty open subset of $X$.

By $k$-type mixing of $(X,T)$, for each pair $(U_i,V_i)$ there exists a finite set
$F_i \subset \{m \in \mathbb{Z}^d : m \ge^k 0\}$ such that
\[
T^{m}(U_i) \cap V_i \neq \emptyset
\quad \text{for all } m \ge^k 0 \text{ with } m \notin F_i.
\]
Let
\[
F := \bigcup_{i=1}^{n} F_i,
\]
which is a finite subset of $\{m \in \mathbb{Z}^d : m \ge^k 0\}$.

For any $m \ge^k 0$ with $m \notin F$, we have
$T^{m}(U_i) \cap V_i \neq \emptyset$ for all $i=1,\dots,n$.
Hence each set $U_i \cap T^{-m}(V_i)$ is non-empty and open.
Define
\[
G := \VOpen{U_1 \cap T^{-m}(V_1),\, U_2 \cap T^{-m}(V_2),\, \dots,\, U_n \cap T^{-m}(V_n)}.
\]
Then $G$ is a non-empty open subset of $\SubX$.

If $A \in G$, then
\begin{enumerate}
    \item $A \cap U_i \neq \emptyset$ for all $i$, hence $A \in \mathcal{U}$, and
    \item $A \cap T^{-m}(V_i) \neq \emptyset$ for all $i$, which implies
    $T^{m}(A) \cap V_i \neq \emptyset$ and therefore $\That^{m}(A) \in \mathcal{V}$.
\end{enumerate}
Thus,
\[
\That^{m}(\mathcal{U}) \cap \mathcal{V} \neq \emptyset
\quad \text{for all } m \ge^k 0 \text{ with } m \notin F,
\]
showing that $(\SubX,\That)$ is $k$-type mixing.

\medskip

Conversely, assume that $(\SubX,\That)$ is $k$-type mixing.
Let $U,V \subseteq X$ be non-empty open sets and define
\[
\mathcal{U} := \VOpen{U},
\qquad
\mathcal{V} := \VOpen{V}.
\]
Then $\mathcal{U}$ and $\mathcal{V}$ are non-empty open subsets of $\SubX$.
By $k$-type mixing of $(\SubX,\That)$, there exists a finite set
$F \subset \{m \in \mathbb{Z}^d : m \ge^k 0\}$ such that
\[
\That^{m}(\mathcal{U}) \cap \mathcal{V} \neq \emptyset
\quad \text{for all } m \ge^k 0 \text{ with } m \notin F.
\]
For such an $m$, there exists $A \in \mathcal{U}$ with $\That^{m}(A) \in \mathcal{V}$.
By definition, $A \subseteq U$ and $T^{m}(A) \subseteq V$, which implies
\[
T^{m}(U) \cap V \neq \emptyset.
\]
Since this holds for all $m \ge^k 0$ with $m \notin F$, we conclude that $(X,T)$ is $k$-type mixing.
\end{proof}

\subsection{$k$-type Li-Yorke Chaos} We first establish few results about $k$-type proximal pairs and $k$-type asymptotic pairs and then we discuss $k$-type Li-Yorke chaos.

\begin{theorem}
\label{thm:k_prox_equivalence}
Let $(X,T)$ be a $\mathbb{Z}^d$-action on a compact metric space $X$,
and let $(\SubX,\That)$ be the induced action on the hyperspace
$\SubX$.

\begin{enumerate}
\item If $(x,y)\in k\text{-Prox}(T)$, then
\(
(\{x\},\{y\})\in k\text{-Prox}(\That).
\)

\item If $(A,B)\in k\text{-Prox}(\That)$, then there exist
$x\in A$ and $y\in B$ such that
\(
(x,y)\in k\text{-Prox}(T).
\)
\end{enumerate}
\end{theorem}

\begin{proof}
Suppose $(x,y)\in k\text{-Prox}(T)$.
Then there exists a sequence $(t_s)$ such that $t_{s+1}>^kt_{s}$ and  
\[
d(T^{t_s}(x),T^{t_s}(y))\to 0.
\]
Since
\[
H(T^{t_s}(\{x\}),T^{t_s}(\{y\}))
=
d(T^{t_s}(x),T^{t_s}(y)),
\]
it follows that
\[
H(T^{t_s}(\{x\}),T^{t_s}(\{y\}))\to 0.
\]
Hence $(\{x\},\{y\})\in k\text{-Prox}(\That)$.

\medskip

Suppose $(A,B)\in k\text{-Prox}(\That)$.
Then there exists a sequence $t_{s+1}>^kt_{s}$ such that
\[
H(T^{t_s}(A),T^{t_s}(B))\to 0.
\]
For each $s$, choose $a_s\in A$ and $b_s\in B$
such that
\[
d(T^{t_s}(a_s),T^{t_s}(b_s))
\le
H(T^{t_s}(A),T^{t_s}(B)).
\]
Since $A$ and $B$ are compact,
we may pass to subsequences (still denoted by $a_s$ and $b_s$)
such that
\[
a_s\to a\in A,
\qquad
b_s\to b\in B.
\]
By continuity of $T$,
\[
d(T^{t_s}(a),T^{t_s}(b))\to 0.
\]
Therefore $(a,b)\in k\text{-Prox}(T)$.
\end{proof}

\begin{theorem} \label{thm:k_asym_equivalence} Let $(X,T)$ be a $\mathbb{Z}^d$-action on a compact metric space $X$, and let $(\SubX,\That)$ be the induced action on $\SubX$. Then,  $$(x,y)\in k\text{-Asym}(T) \text{ if and only if }  (\{x\},\{y\})\in k\text{-Asym}(\That).$$  \end{theorem}
\begin{proof}

Suppose $(x,y)\in k\text{-Asym}(T)$.
By definition, for every $\epsilon>0$ there exists $n>^k0$ such that
\[
d(T^{n+m}(x),T^{n+m}(y))\le \epsilon
\quad \text{for all } m>^k0.
\]
Since the Hausdorff distance between singleton sets coincides
with the metric, we have
\[
H\bigl(T^{n+m}(\{x\}),T^{n+m}(\{y\})\bigr)
=
d(T^{n+m}(x),T^{n+m}(y))\le \epsilon
\]
for all $m>^k0$.
Hence
\[
(\{x\},\{y\})\in k\text{-Asym}(\That).
\]

\medskip

Suppose $(\{x\},\{y\})\in k\text{-Asym}(\That)$.
Then for every $\epsilon>0$ there exists $n>^k0$ such that
\[
H\bigl(T^{n+m}(\{x\}),T^{n+m}(\{y\})\bigr)\le \epsilon
\quad \text{for all } m>^k0.
\]
Again using the equality of Hausdorff and pointwise distance for
singletons, we obtain
\[
d(T^{n+m}(x),T^{n+m}(y))\le \epsilon
\quad \text{for all } m>^k0.
\]
Thus $(x,y)\in k\text{-Asym}(T)$.
\end{proof}

\begin{theorem}
\label{thm:k_LY_chaos_equivalence}
Let $(X,T)$ be a $\mathbb{Z}^d$-action on a compact metric space $X$,
and let $(\SubX,\That)$ be the induced action on the hyperspace
$\SubX$.
If $(X,T)$ is $k$-type Li--Yorke chaotic then
$(\SubX,\That)$ is $k$-type Li--Yorke chaotic.
\end{theorem}

\begin{proof}\noindent Assume that $(X,T)$ is $k$-type Li--Yorke chaotic.
Let $S\subset X$ be an uncountable $k$-scrambled set.
Consider the collection
\[
\mathcal{S} := \bigl\{ \{x\} : x\in S \bigr\}
\subset \SubX.
\]
Since the map $x\mapsto \{x\}$ is injective,
$\mathcal{S}$ is uncountable.

Let $\{x\},\{y\}\in\mathcal{S}$ with $x\neq y$.
Then $(x,y)\in k\text{-Prox}(T)$ implies that $ (\{x\},\{y\})\in k\text{-Prox}(\That)$ by Theorem \ref{thm:k_prox_equivalence}.

If $(\{x\},\{y\})$ were in $k\text{-Asym}(\That)$,
then by Theorem \ref{thm:k_asym_equivalence}, it follows that
$(x,y)\in k\text{-Asym}(T)$,
a contradiction.
Hence
$(\{x\},\{y\})
\notin k\text{-Asym}(\That).$

Thus every distinct pair in $\mathcal{S}$
forms a $k$-type Li--Yorke pair for $(\SubX,\That)$,
so $\mathcal{S}$ is an uncountable $k$-scrambled set
in $\SubX$.
Therefore $(\SubX,\That)$ is $k$-type Li--Yorke chaotic.
\end{proof}

\noindent\textbf{Acknowledgements.} {The author gratefully acknowledges the financial support from the Council of Scientific and Industrial Research (CSIR), Government of India, through the fellowship file no. 09/1026(0044)/2021-EMR-I. The author also acknowledges the guidance of Prof. Sharan Gopal in framing this article.}

\end{document}